\documentclass[10pt,reqno]{amsart}
\usepackage{amssymb,accents}
\def\noi{\noindent}

\newtheorem{Thm}{Theorem}[section]

\newtheorem{Lm}[Thm]{Lemma}

\newtheorem{Cor}[Thm]{Corollary}

\newtheorem{state}{Theorem}

\def\cal{\mathcal}
\def\Bbb{\mathbb}
\def\mf{\mathfrak}

\def\<{\langle}
\def\>{\rangle}

\def\a{\alpha}

\def\d{\delta}

\def\th{\theta}

\def\l{\lambda}
\def\L{\Lambda}

\def\Re{\Bbb R}
\def\F{\Bbb F}
\def\C{\Bbb C}
\def\Z{\Bbb Z}

\def\W{\ring W}
\def\w{\ring w}
\def\RR{\ring R}
\def\Q{\ring Q}

\begin{document}
\title[]
{A weight multiplicity formula for Demazure modules}
\author{Bogdan Ion}\thanks{Supported in part by a Rackham Faculty
Research Fellowship}
\address{Department of Mathematics, University of Michigan, Ann
Arbor, MI 48109} \email{bogdion@umich.edu}
\begin{abstract}
We establish a formula for the weight multiplicities of Demazure
modules (in particular for highest weight representations) of a
complex connected algebraic group in terms of the geometry of its
Langlands dual.
\end{abstract}
\maketitle

\thispagestyle{empty}
\section*{Introduction}

Let $G$ be a complex, connected, simple algebraic group  and let
$B$ and $T$ be a Borel subgroup and maximal torus of $G$ contained
in $B$. For an integral dominant weight $\l$, denote by $V_\l$ the
integrable highest weight $G$-module with highest weight $\l$. The
Demazure modules $D_{w(\l)}$ (with $w$ running over the Weyl group
of $G$) form a filtration of $V_\l$. The module $D_{w(\l)}$ is
defined as the $B$--submodule of $V_\l$ generated by the extremal
weight space of weight $w(\l)$ (which are one--dimensional). In
particular, if $w_\circ$ is the longest element in the Weyl group
of $G$, $D_{w_\circ(\l)}$ is precisely $V_\l$. More importantly,
the Demazure modules can be realized as spaces of global sections
of certain line bundles over Schubert varieties (see Section
\ref{demazure} for details).

The main goal of this note is to explain how the recent
developments \cite{ion2}, \cite{ion3} in the theory of
nonsymmetric Macdonald polynomials have as consequence a geometric
formula for weight multiplicities in Demazure modules and, in
particular, irreducible $G$--modules. We will briefly state our
result.

Let ${G^\vee}$ be the unique complex, connected, reductive group
whose root datum is dual to that of $G$.  Let us also consider
$B^\vee$, a Borel subgroup of $G^\vee$. The group $G^\vee$ is also
defined over $\cal O:=\C[[x]]$. We denote $G^\vee(\cal O)$ by $K$
and let $I$ be the subgroup of $K$ defined as the inverse image of
$B^\vee(\C)$ under the reduction map $K\to G^\vee(\C)$. The space
$G^\vee(\C((x)))/I$ is endowed with a structure of ind--variety
over $\C$. For any integral weights $\l$ and $\mu$ let
$M_{\l,\mu}(\C)$ be the finite dimensional subvariety of
$G^\vee(\C((x)))/I$ defined in Section \ref{varieties}.

\begin{state}\label{thm1}
Let $\l$ and $\mu$ be two integral weights such that $\mu\leq \l$.
The weight multiplicity $m_{\l,\mu}$ of the weight $\l$ in the
Demazure module $D_\l$ equals the number of top dimensional
irreducible components of $M_{\l,\mu}(\C)$.
\end{state}

Other geometric formulas for weight multiplicities in irreducible
$G$--modules exist in the literature \cite{lusztig}, \cite{mv} but
no such formula was available for Demazure modules. Our result is
similar to one consequence \cite[Corollary 7.4]{mv} of the recent
work of Mirkovi\' c and Vilonen  on the geometric Satake
isomorphism (see \ref{remarks} on further comments about this
connection). However, we should keep in mind that the Demazure
modules are generally just $B$--modules and their result covers
only the case of reductive groups; it would be interesting to
investigate if the ideas in \cite{mv} could be used to cover a
larger class of groups.

\noi {\sl Acknowledgement:} I would like to thank Brian Conrad for
a clarifying conversation.
\section{Preliminaries}\label{prelim}

\subsection{}
Let ${\mf g}$ be the complex, simple Lie algebra of rank $n$ and
let $\mf h$ and $\mf b$ be a Cartan subalgebra and, respectively,
a Borel subalgebra containing $\mf h$. Let $\RR$ be root system of
$\mf g$ with respect to $\mf h$ and denote by $R^+$ the set of
roots of $\mf b$ with respect to $\mf h$. The set of positive
simple roots is denoted by $\{\a_1,\dots,\a_n \}$ and
$\{\l_1,\cdots,\l_n\}$ will denote the fundamental weights.  As
usual, $\Q$ denotes the root lattice of $\RR$ and $P$ denotes its
weight lattice.

The vector space ${\mf h}_\Re^*$ (the real vector space spanned by
the roots) has a canonical scalar product $(\cdot, \cdot)$ which
we normalize such that it gives square length 2 to the short roots
in $\RR$ (if there is only one root length we consider all roots
to be short). We will use $\RR_s$ and $\RR_\ell$ to refer to the
short and respectively long roots in $\RR$. We will identify the
vector space $\mf h_\Re$ (the real vector space spanned by the
coroots) with its dual using this scalar product. Under this
identification $\a^\vee=2\a/(\a,\a)$ for any root $\a$. The root
$\th$ is defined as the {\sl highest short root} in $\RR$. Also,
let us consider $$\rho=\frac{1}{2}\sum_{\a\in R}\a^\vee$$

To the finite root system $\RR$ we will associate an {\sl affine
root system} $R$. Let ${\rm Aff}(\mf h_\Re)$ be the space of
affine linear transformations of $\mf h_\Re$. As a vector space,
it can be identified to $\mf h_\Re^*\oplus \Re\d$ via
$$
(f+c\d)(x)=f(x)+c, \quad \text{for~} f\in \mf h_\Re^*, ~x\in\mf
h_\Re \text{~and~} c\in \Re
$$
Let $r$ denote the maximal number of laces connecting two vertices
in the Dynkin diagram of $\RR$. Then,
$$R:=(\RR_s+\Z\d)\cup(\RR_\ell +r\Z\d)\subset \mf h_\Re^*\oplus \Re\d$$

The set of affine positive roots $R^+$ consists of affine roots of
the form $\a+k\d$ such that $k$ is non--negative if $\a$ is a
positive root, and $k$ is strictly positive if $\a$ is a negative
root. The affine simple roots are $\{\a_i\}_{0\leq i\leq n}$,
where we set $\a_0:=\d-\th$.


\subsection{}
The scalar product on $\mf h_\Re^*$ can be extended to a
non--degenerate bilinear form on the real vector space
$$H:=\mf h_\Re^*\oplus \Re\d\oplus \Re\L_0$$ by requiring that $(\d,\mf h_\Re^*\oplus \Re\d)=
(\L_0,\mf h_\Re^*\oplus \Re\L_0)=0$ and $(\d,\L_0)=1$. Given
$\a\in R$ and $x\in H$ let
$$
s_\a(x):=x-\frac{2(x,\a)}{(\a,\a)}\a\
$$
The {\sl affine Weyl group} $W$ is the subgroup of ${\rm GL}(H)$
generated by all $s_\a$ (the simple reflections $s_i=s_{\a_i}$ are
enough). The {\sl finite Weyl group} $\W$ is the subgroup
generated by $s_1,\dots,s_n$. It is easy to see that $\mf
h_\Re^*+\Re\d+\L_0$ is stable under the action of $W$. Therefore,
if we identify $\mf h_\Re^*$ with $\mf h_\Re^*+\Re\d+\L_0$ we
obtain an affine action of $W$ on $\mf h_\Re^*$;  we denote by
$w\cdot x$ the affine action of $w\in W$ on  $x\in \mf h_\Re^*$.
For example, the affine  action of $s_0$ can be described as
follows
\begin{equation}
s_0\cdot x   =s_\th(x)+\th\
\end{equation}

We define the fundamental affine chamber as
\begin{equation}\label{affinechamber}
\cal C:=\{ x\in \mf h_\Re^*\ | \ (x+\L_0,\a_i^\vee)\geq 0\ ,\
0\leq i\leq n\}
\end{equation}
The non-zero elements of $\cal O_P:=P\bigcap \cal C$ are the
so-called minuscule weights. Let us remark that each orbit of the
affine action of $W$ on $P$ contains the origin or a unique
element of $\cal O_P$.

\subsection{}

For each $w$ in $W$ let $\ell(w)$ be the length of a reduced (i.e.
shortest) decomposition of $w$ in terms of simple reflections. The
length of $w$ can be also geometrically described as follows. For
any affine root $\a$, denote by $H_\a$ the affine hyperplane
consisting of fixed points of the affine action of $s_\a$ on $\mf
h_\Re^*$. Then, $\ell(w)$ equals the number of affine hyperplanes
$H_\a$ separating $\cal C$ and $w\cdot\cal C$. For any affine
transformation of $\mf h_\Re^*$ which preserves the set of
hyperplanes $\{H_\a\}_{\a\in R}$, we can use the geometric point
of view to define the length of that transformation. For example,
for a weight $\l$ we define the following affine transformation of
$\mf h_\Re^*$
$$
\tau_\l(x)=x+\l
$$
It is easy to check that $\tau_\l$ has the required properties to
allow us to talk about its length.  In fact a concrete formula for
its length is available (see, for example, \cite[(5)]{ion2})
\begin{equation}\label{length}
\ell(\tau_\l)=\sum_{\a\in\RR} |(\l,\a^\vee)|
\end{equation}

For each weight $\l$ define $\l_-$, respectively $\tilde\l$, to be
the unique element in $\W(\l)$, respectively $W\cdot\l$, which is
an anti-dominant weight, respectively an element of $\cal O_P$
(that is a minuscule weight or zero), and
$\w_\l\hspace{-0.2cm}^{-1}\in\ \W$, respectively $w_\l^{-1}\in W$,
to be the unique minimal length element by which this is achieved.
Also, for each weight $\l$ define $\l_+$ to be the unique element
in $\W(\l)$ which is dominant and denote by $w_\circ$ the maximal
length element in $\W$. It was shown in \cite[Lemma 1.7 (3)]{ion2}
that the following equality holds for any weight $\l$
\begin{equation}\label{eq3}
\ell(\tau_\l)=\ell(w_\l)+\ell(\w_\l)
\end{equation}
For later use we record the following.
\begin{Lm}\label{lemma1}
Let $\l$ be a weight. Then
$$
2\<\l,\rho\>\leq \ell(w_\l)+\ell(\w_\l)
$$
\end{Lm}
\begin{proof}
From (\ref{length}) we know that $2\<\l,\rho\>\leq \ell(\tau_\l)$.
Now, (\ref{eq3}) implies the desired result.
\end{proof}


\subsection{}
The Bruhat order is a partial order on any Coxeter group defined
in way compatible with the length function. For an element $w$ we
put $w<s_iw$\  if and only if \ $\ell(w)<\ell(s_iw)$. The
transitive closure of this relation is called the Bruhat order.
The terminology is motivated by the way this ordering arises for
Weyl groups in connection with inclusions among closures of Bruhat
cells for a corresponding semisimple algebraic group.

We can use the Bruhat order on $W$ to define a partial order on
the weight lattice which will also be called the Bruhat order. For
any $\l,\mu\in P$ we write
\begin{equation}\label{order}
\l<\mu\ \ \text{if and only if}\ \ w_\l<w_\mu
\end{equation}
The minimal elements of the weight lattice with respect to this
partial order are the minuscule weights and if $\l<\mu$ then
necessarily $\tilde \l=\tilde \mu$.

\begin{Lm}\label{lemma2}
Let $\l$ and $\mu$ be two weights such that $\mu\leq \l$. Define
the rational number
\begin{equation}\label{n}
n_{\l,\mu}:=\frac{1}{2}\ell(w_\l)-\frac{1}{2}\ell(\w_\l)+\ell(w_\circ)-\<\mu,\rho\>
\end{equation}
Then $n_{\l,\mu}$ is a positive integer.
\end{Lm}
\begin{proof} Let us argue first that $n_{\l,\mu}$ is positive.
We can use Lemma \ref{lemma1} to obtain
$$
n_{\l,\mu}\geq
\frac{1}{2}\ell(w_\l)-\frac{1}{2}\ell(\w_\l)+\ell(w_\circ)-
\frac{1}{2}\ell(w_\mu)-\frac{1}{2}\ell(\w_\mu)
$$
Keeping in mind that $\ell(w_\l)\geq\ell(w_\mu)$ (which is a
consequence of the hypothesis) and that $w_\circ$ is the maximal
length element in $\W$, our claim immediately follows.

To show that $n_{\l,\mu}$ is integer, let us remark that it is
enough to check that
$\frac{1}{2}\ell(w_\l)+\frac{1}{2}\ell(\w_\l)$ is integer. This
fact follows from (\ref{eq3}) together with
$$\ell(\tau_\l)=\ell(\tau_{\l_+})=2\<\l_+,\rho\>$$
\end{proof}

\section{Demazure modules}\label{demazure}

Let $G$ be a complex, connected, simple algebraic group with Lie
algebra $\mf g$ and denote by $T$ and $B$ the maximal torus  and
the Borel subgroup of $G$ with Lie algebras $\mf h$ and $\mf b$,
respectively.

For an integral dominant weight $\l$, denote by $V_\l$ the
integrable highest weight ${\mf g}$-module with highest weight
$\l$. It is well known that $V_\l$ is an irreducible $\mf
g$--module. Hence, for any $w$ in $\W$ the $w(\l)$--weight space
$V_{\l,w(\l)}$ is one--dimensional. The {\sl Demazure module}
$D_{w(\l)}$ is defined as the ${\mf b}$--module generated by
$V_{\l, w(\l)}$. Since $\l$ is integral $D_{w(\l)}$ is also a
$B$--module. The Demazure modules associated to a fixed integral
dominant weight $\l$ form a filtration (with respect to the Bruhat
order on $\W$) of $V_\l$. In particular, $D_\l=V_{\l,\l}$ and
$D_{w_\circ(\l)}=V_\l$.

There exists an important geometrical construction of Demazure
modules which relates them with the geometry of Schubert varieties
in the  flag variety $G/B$. For any $w$ in $\W$ let
$$
X_w=\overline{BwB/B}\subseteq G/B
$$
denote the Schubert variety associated to $w$. The Schubert
varieties are closed finite dimensional projective irreducible
subvarieties of $G/B$. For any  integral weight $\l$ we denote by
$e^\l$ the character of $B$ obtained from $T$ via the isomorphism
$T\simeq B/[B,B]$. Consider the fiber product $\cal
L_\l=G\times_B\C_\l$, where $\C_\l$ denotes $\C$ equipped with the
$B$--action given by the character $e^\l$. The natural projection
$G\times\C_\l\to G$ induces a well defined $G$--equivariant
holomorphic map $\cal L_\l\to G/B$; in other words $\cal L_\l$
becomes a $G$--equivariant holomorphic line bundle over $G/B$. By
restriction we obtain a line bundle $\cal L_{\l,w}$ over $X_w$.
Since $X_w$ is $B$--invariant, the space of holomorphic sections
$H^0(X_w,\cal L_{\l,w})$ admits a $B$--module structure. The
relation with the Demazure module $D_{w(\l)}$ is the following.
\begin{Thm} Let $\l$ be a dominant integral weight and let $w$ be an element of $W$.
Then, $D_{w(\l)}$ and $H^0(X_w,\cal L_{-\l, w})^*$ are isomorphic
as $B$--modules.
\end{Thm}
The result holds in the more general setting of generalized flag
varieties of Kac--Moody groups.  For a proof see, for example,
Corollary 8.1.26 in \cite{kumar}.

Let $\l$ be an arbitrary integral weight. As $T\subset B$, the
Demazure module $D_{\l}$  is also a $T$--module. Its
$T$--character will be denoted by $\chi_\l$.  Let us write
\begin{equation}\label{eq1}
\chi_\l=\sum_{\mu\in X^*(T)}m_{\l,\mu}e^\mu
\end{equation}
where we denoted by $X^*(T)$ the character group of the torus $T$.
The non--negative integers $m_{\l,\mu}$ appearing in the above
formula are the multiplicities of the weights $\mu$ in the
Demazure module $D_\l$.  As remarked before, if $\l$ is dominant,
$\chi_\l$ equals $e^\l$ and $\chi_{w_\circ(\l)}$ is the character
of $V_\l$.

\section{Nonsymmetric Macdonald polynomials}\label{nonsymm}

\subsection{}

Let us introduce a field $\C_{q,t}$ (of parameters) as follows.
Let $q$ and $t$ be two formal parameters and let $m$ be the lowest
common denominator of the rational numbers $\{(\a_j,\l_k)\ |\
1\leq j,k\leq n \}$. The field $\C_{q,t}$ is defined as the field
of rational functions in $q^{\frac{1}{m}}$ and $t^{\frac{1}{2}}$
with complex coefficients. We will also use the field of rational
functions in $t^{\frac{1}{2}}$ denoted by $\C_t$. The algebra
$\cal R_{q,t}=\C_{q,t}[e^\l;\l\in P]$ is the group
$\C_{q,t}$-algebra of the lattice $P$. Similarly, the algebra
$\cal R_t=\C_t[e^\l;\l\in P]$ is the group $\C_t$-algebra of the
lattice $P$.

The nonsymmetric Macdonald polynomials $E_\l(q,t)$ associated to
the root system $\RR$  are remarkable family of  polynomials
indexed by the weight lattice $P$ and which form a basis for $\cal
R_{q,t}$. They were defined by Opdam and Macdonald for various
specializations of the parameters and in full generality by
Cherednik \cite{c}. We refer the reader to \cite{c} for a detailed
account of their construction and basic properties. We only
mention at this point that they satisfy the following
triangularity property with respect to the Bruhat order
$$
E_\l(q,t)\in {\rm span}\<~e^\mu~|~\mu\leq\l~\>
$$


\subsection{}
Although, a priori, the coefficients of the polynomials
$E_\l(q,t)$ are just rational functions in $q$ and $t$  it was
proved in \cite[Section 3.1]{ion2} that in fact each coefficient
has a finite limit as $q$ approaches infinity and, in consequence,
$$E_\l(t)=\lim_{q\to \infty}E_\l(q,t)$$ is well defined as an
element of $\cal R_t$.

Another important fact is that the coefficients of $E_\l(t)$ are
polynomials in $t^{-1}$ and in consequence their limit as $t$
approaches infinity is finite. Moreover, the following is true.
Let $G$, $B$ and $T$ be a complex connected simple algebraic
group, a Borel subgroup and a maximal torus as in Section
\ref{demazure}. We can certainly regard $X^*(T)$, the character
group of the torus $T$, as being a sublattice of $P$.  Therefore,
the Demazure character $\chi_\l$ associated to an integral weight
$\l$ becomes an element of $\cal R_t$, and the formula (\ref{eq1})
will be regarded as an identity in $\cal R_t$. The following
result was proved in \cite[Corollary 3.8]{ion2}.

\begin{Thm}\label{limit} Let $\l$ be an arbitrary integral weight for $G$.
Then,
\begin{equation}
\chi_\l=\lim_{t\to\infty}E_\l(t)
\end{equation}
\end{Thm}


\subsection{} For special values of the parameter $t$ the
coefficients of the polynomials $E_\l(t)$ have a rather different
interpretation. To be able to state the result we first need to
introduce more notation.

Given  $\l$ in $P$ let us define the following normalization
factor
\begin{equation}\label{eq4}
j_\l=
t^{\frac{1}{2}(\ell(w_\l)-\ell(\w_\l))}
\end{equation}
where $w_\l$ and $\w_\l$ are the Weyl group elements defined in
Section \ref{prelim}. Let us remark that, as it follows from
\cite[Corollary 3.4]{ion2}, the normalization factor denoted by
the same symbol in \cite[(12)]{ion3} is, for equal values of the
parameters, precisely the element defined above.

\subsection{}

Let ${\bf G^\vee}$ be the unique split, connected, reductive group
scheme whose root datum is dual to that of $G$ (the Chevalley
group scheme).  Let ${\bf B^\vee}$ be a Borel subgroup of ${\bf
G^\vee}$ and ${\bf T^\vee}$ a maximal split torus of $\bf G^\vee$
contained in ${\bf B^\vee}$. The Borel ${\bf B^\vee}$ has the Levi
decomposition ${\bf B=T^\vee U}$, where ${\bf U}$ is the unipotent
radical of ${\bf B^\vee}$. The unique Borel subgroup of $\bf
G^\vee$ which is opposed to $\bf B^\vee$ with respect to $\bf
T^\vee$ will be denoted by $\overline{\bf B}^\vee$ and
$\overline{\bf U}$ denotes its unipotent radical. We have
deliberately ignored the field of definition since all above
groups are defined over $\Z$ (the structure constants involve only
integers).

Let $p$ be a prime number and $\mf t$ a positive integer power of
$p$. We will denote by $\F_{\mf t}$ the finite field with $\mf t$
elements. For the moment let $k$  denote an arbitrary field. Let
$x$ be a formal parameter and let $F:=k((x))$ be the quotient
field of $\cal O:=k[[x]]$ (formal power series with coefficients
in $k$). Of course, $F$ is a $p$--adic field, $\cal O$ is its ring
of integers and $k$ is the residue field.

The $F$--rational points of ${\bf G^\vee}$ will be denoted by
$G^\vee(F)$ (or simply by $G^\vee$ if the field $F$ is implicitly
understood) and the same type of notation will be used for all the
linear algebraic groups defined above.  Each $\l\in X^*(T)$
becomes an element of $X_*(T^\vee)$ (the cocharacter group of
$T^\vee$) and therefore determines a morphism $\l:F^\times\to
T^\vee$. We will denote by $x^\l$ the image of $x$ under the above
morphism. The groups $\bf G^\vee$, $\bf T^\vee$, $\bf U$,
$\overline{\bf U} $ are also defined over $\cal O$. We will denote
$G^\vee(\cal O)$ by $K$ and denote by $I$ the Iwahori subgroup
defined as the inverse image of $B^\vee(k)$ under $G^\vee(\cal
O)\to G^\vee(k)$. Of course, all the above considerations hold if
we replace $k$ by $\Z$.

\subsection{}
In this subsection we will assume that the field $k$ is $\F_{\mf
t}$. Let us remark that the affine root system associated to
$(G^\vee,T^\vee)$ is in fact $R^\vee$ and that $K$ is a maximal
compact subgroup of $G^\vee$.

The group $G^\vee$ is unimodular as opposed to the Borel subgroup
$B^\vee$ which is not. We choose a Haar measure on $G^\vee$
normalized such that the Iwahori subgroup $I$ has volume one. The
modular function of the Borel subgroup $\d:T^\vee\to \Re^\times_+$
is defined by
 $\d(a)=|{\rm det}(Ad_{|U}(a))|$ for any $a$ in $T^\vee$;
 we denoted by $Ad_{|U}$ the automorphism of the Lie algebra of
$U$ given by the adjoint representation and we denoted by
$|\cdot|$ the usual metric on $F$ induced by the valuation. Since
$G^\vee$ is split, the formula for the modular function takes the
following form on the elements $x^\l$, for $\l\in X^*(T)$
\begin{equation}\label{eq2}
\d{}(x^\l)={\mf t}^{2\<\l,\rho\>}
\end{equation}

For an element $f=f(t)$ in $\cal R_t$ we will write $f(\mf t)$ to
refer to the element of $\C[e^\l;\l\in P]$ obtained by
substituting the positive integer $\mf t$ for the parameter $t$.
We are now ready to state a result connecting the coefficients of
nonsymmetric Macdonald polynomials with the geometry of the group
$G^\vee$. The  Theorem stated below was proved in \cite{ion3} (see
Theorem 5.10 and formula (25)).
\begin{Thm}\label{padic}
Let $\l$ be an element of $X^*(T)$ (or, equivalently, of
$X_*(T^\vee))$. The coefficients appearing in
$$
E_\l(\mf t)=\sum_{\mu\leq \l} c_{\l,\mu} e^\mu
$$
are given by
$$
c_{\l,\mu}=\frac{vol(\overline{U} x^{-\mu}I\cap
Kx^{-\l}I)}{j_\l(\mf t)\mf
t^{\ell(w_\circ)}\d^{\frac{1}{2}}(x^{w_\circ(\mu)})}
$$
\end{Thm}
We would prefer to make the denominator in the above formula as
explicit as possible.
\begin{Lm}
Let $\l$ and $\mu$ be two elements of $X^*(T)$ for which $\mu\leq
\l$. Then,
$$
j_\l(\mf t)\mf
t^{\ell(w_\circ)}\d^{\frac{1}{2}}(x^{w_\circ(\mu)})=\mf
t^{n_{\l,\mu}}
$$
\end{Lm}
\begin{proof}
Straightforward from (\ref{n}), (\ref{eq4}) and (\ref{eq2}).
\end{proof}

 By combining Theorem \ref{limit} and Theorem \ref{padic} we can see that the
 weight multiplicities (as defined by (\ref{eq1})) of the Demazure module
 $D_\l$ can be computed as follows.
 \begin{Cor}\label{cor1}
Let $\l$ and $\mu$ be two elements of $X^*(T)$ for which $\mu\leq
\l$. Then,
 \begin{equation}\label{eq5}
m_{\l,\mu}=\lim_{\mf t\to\infty}\frac{vol(\overline{U}
x^{-\mu}I\cap Kx^{-\l}I)}{\mf t^{n_{\l,\mu}}}
\end{equation}
\end{Cor}
\section{The multiplicity formula}
\subsection{The varieties $ M_{\l,\mu}$}\label{varieties}

Let us assume for a moment that $k=\C$. It is well--known (see,
for example, \cite{mv} and the references therein) that the space
$G^\vee/K$ is an {\sl ind--variety} (i.e. admits an increasing
filtration with varieties such all inclusion maps among them are
closed embeddings) defined over $\C$. We refer the reader to
\cite[Chapter IV]{kumar} for a brief introduction to
ind--varieties. In our case, the members of the filtration on
$G^\vee/K$ can be constructed from a filtration of $G^\vee$
obtained by bounding the number of poles with respect to $x$ of
the matrix coefficients for a faithful representation of $G^\vee$.

The space $ G^\vee/K$ is usually referred to as the {\sl affine
Grassmannian} of $G^\vee(\C)$. However, because $ G^\vee$ and $K$
are still defined for $k=\Z$, the same geometric constructions go
through in this case, and $G^\vee/K$ acquires a structure of
ind--variety defined over $\Z$ (see also the remarks in Section 14
of \cite{mv}). In a completely similar way, keeping in mind that
$I$ is defined for $k=\Z$, we endow $G^\vee/I$ with a structure of
ind--variety defined over $\Z$.

Let us consider the map $\pi:G^\vee\to G^\vee/I$. The group $K$
acts on $G^\vee/I$ by finite dimensional orbits. Indeed, the
orbits are all of the form $\pi(Kx^\nu I)$, and the number of
poles of the elements in $Kx^\nu I$ is bounded by $$\max_{1\leq
i\leq n}|\<\nu,\a_i^\vee\>|$$ Therefore, $\pi(Kx^\nu I)$ is
included in one member of the filtration on $G^\vee/I$ and hence
it is finite dimensional. For any $\l$ and $\mu$ in $X_*(T)$, let
us define the variety $M_{\l,\mu}$ as $\pi(\overline{U}
x^{-\mu}I\cap Kx^{-\l}I)$. From the above considerations it is
clear that the varieties $M_{\l,\mu}$ are finite dimensional and
defined over $\Z$.


\subsection{Proof of Theorem \ref{thm1}}
We will show that the right hand side of the formula (\ref{eq5})
equals the number of irreducible components of top dimension of
$M_{\l,\mu}(\C)$. Let us recall that Corollary \ref{cor1} holds
under the hypothesis $k=\F_{\mf t}$. Because the Iwahori subgroup
$I$ has volume one we can regard the volume of the set
$\overline{U} x^{-\mu}I\cap Kx^{-\l}I$ as the number of right
$I$--cosets in $\overline{U} x^{-\mu}I\cap Kx^{-\l}I$ or,
equivalently, as the number of points in $M_{\l,\mu}(\F_{\mf t})$
(which we denote by $|M_{\l,\mu}(\F_{\mf t})|$).

The Lefschetz fixed point formula (for the Frobenius automorphism)
and Theorem 1 in \cite{deligne} (for $\cal F_0=\Bbb Q_\ell$,
$S_0=Spec(\F_{\mf t})$ and $X_0=M_{\l,\mu}$)  let us conclude that
$$\lim_{\mf t\to\infty}\frac{|M_{\l,\mu}(\F_{\mf t})|}{\mf
t^{n_{\l,\mu}}}$$  is indeed the number of irreducible components
of top dimension of $M_{\l,\mu}(\C)$.

As an immediate consequence we obtain the following.
\begin{Cor}
Let $\l$ and $\mu$ be two elements of $X^*(T)$ for which $\mu\leq
\l$. The dimension of $M_{\l,\mu}(\C)$ is $n_{\l,\mu}$.
\end{Cor}

\subsection{Final remarks}\label{remarks}

If we keep in mind that for anti--dominant $\l$ the Demazure
module $D_\l$ is in fact $V_{\l_+}$, the irreducible $G$--module
with highest weight $\l_+=w_\circ(\l)$ (or lowest weight $\l$),
our result gives a geometric  formula for the weight
multiplicities of $V_{\l_+}$. However, if we would be only
interested in a formula for weight multiplicities in $V_{\l}$
(here $\l$ is dominant), this could be obtained by closely
following the argument given above but replacing in Theorem
\ref{limit} the nonsymmetric Macdonald polynomials with symmetric
Macdonald polynomials $P_\l(t)$ (at $q\to\infty$) and Theorem
\ref{padic} with corresponding result stating the fact that
$P_\l(\mf t)$ arises as the Satake transform of the characteristic
function of the set $Kx^\l K$. If we denote by $Z_{\l,\mu}(\C)$
the $\C$--points of the image of $Ux^\mu K\cap Kx^\l K$ through
the map $$\tilde\pi:G^\vee\to G^\vee/K$$ the corresponding result
read as follows.
\begin{Thm}\label{thm2} Let $\l$ and $\mu$ be two dominant integral weights
such that $\mu\leq \l$. The weight multiplicity $m_{\l,\mu}$ of
the weight $\mu$ in the irreducible highest weight module $V_\l$
equals the number of top dimensional irreducible components of
$Z_{\l,\mu}(\C)$.
\end{Thm}

The above statement is quite close to one of results in \cite{mv}
which we will briefly recall. For the following statements we
assume that $k=\C$. First, Theorem 3.2 in \cite{mv} shows that
$\tilde\pi(Ux^\mu K)\cap \overline{\tilde\pi(Kx^\l K)}$ (the
closure is in $G^\vee/K$) is pure dimensional (which, of course,
implies that $Z_{\l,\mu}$ is pure dimensional). Second, as a
consequence of the equivalence of tensor categories between the
category of representations of $G$ and the category of
$K$--equivariant perverse sheaves on the affine Grassmannian
$G^\vee/K$ with $\C$--coefficients, Corollary 7.4 in \cite{mv}
states that the weight multiplicity $m_{\l,\mu}$ equals the number
of irreducible components of $\tilde\pi(Ux^\mu K)\cap
\overline{\tilde\pi(Kx^\l K)}$. The absence of the closure in
Theorem \ref{thm2} slightly improves this statement.

The main result in \cite{mv} (the above mentioned equivalence of
categories) was proved for reductive groups. Keeping in mind that
Theorem \ref{thm1} holds for Demazure modules (which are just
$B$--modules) and that $B$ is a solvable group, our result seems
to suggest that the results of Mirkovi\' c and Vilonen might
eventually be extended to a larger class of groups and
representations.

\end{document}